\documentclass[10pt,twoside]{article}
\usepackage{amsfonts}
\usepackage{graphicx}
\usepackage{amsmath}
\usepackage{Latex-document}

\def\text{\mbox}

\title{\bf Estimates for the Strong Approximation \vskip -2mm in
 Multidimensional Central \vskip -2mm Limit Theorem\thanks{Research partially
supported by Russian Foundation of Basic Research (RFBR) Grant
02-01--00265, and by RFBR-DFG Grant 99-01--04027.}\vskip 6mm}
\author{A. Yu. Zaitsev\vspace*{-0.5cm}\thanks{St. Petersburg Branch of the Steklov
Mathematical Institute, Fontanka 27, St. Petersburg 191011,
Russia. E-mail: zaitsev@pdmi.ras.ru}}
\date{\vspace{-8mm}}

\begin{document}

\maketitle

\thispagestyle{first} \setcounter{page}{107}

\begin{abstract}
\vskip 3mm In a recent paper the author obtained optimal bounds for the
strong Gaussian approximation of sums of independent $\mathbf{R}^{d}$-valued
random vectors with finite exponential moments. The results may be
considered as generalizations of well-known results of
Koml\'{o}s--Major--Tusn\'{a}dy and Sakhanenko. The dependence of constants
on the dimension $d$ and on distributions of summands is given explicitly.
Some related problems are discussed.

\vskip 4.5mm

\noindent \textbf{2000 Mathematics Subject Classification:} 60F05, 60F15,
60F17.

\noindent \textbf{Keywords and Phrases:} Strong approximation, Prokhorov
distance, Central limit theorem, Sums of independent random vectors.
\end{abstract}

\vskip 12mm

\section{Introduction} \label{section 1}\setzero
\vskip-5mm \hspace{5mm}

Let $X_{1},\dots ,X_{n},\dots $ be mean zero independent $\mathbf{R}^{d}$%
-valued random vectors and $\Bbb{D}_{n}=\mbox{cov}\,S_{n}$ the covariance
operator of the sum $S_{n}=\sum_{i=1}^{n}X_{i}$. By the Central Limit
Theorem, under some simple moment conditions the distribution of normalized
sums $\Bbb{D}_{n}^{-1/2}S_{n}$ is close to the standard Gaussian
distribution. The invariance principle states that, in a sense, the
distribution of the \textit{whole} sequence $\Bbb{D}_{n}^{-1/2}S_{1},\dots ,%
\Bbb{D}_{n}^{-1/2}S_{n},\ldots $ is close to the distribution of the
sequence $\Bbb{D}_{n}^{-1/2}T_{1},\dots $, $\Bbb{D}_{n}^{-1/2}T_{n},\ldots $%
, where $T_{n}=\sum_{i=1}^{n}Y_{i}$ and $Y_{1},\dots ,Y_{n},\ldots $ is a
corresponding sequence of independent Gaussian random vectors (this means
that \thinspace $Y_{i}$ \thinspace has the same mean and the same covariance
operator as \thinspace $X_{i}$, \thinspace $i=1,\dots ,n,\ldots $).

We consider here the problem of strong approximation which is more delicate
than that of estimating the closeness of distributions. It is required to
construct on a probability space a sequence of independent random vectors
\thinspace $X_{1},\dots ,X_{n}$ \thinspace (with given distributions) and a
corresponding sequence of independent Gaussian random vectors $Y_{1},\dots
,Y_{n}$ so that the quantity
\[
\Delta (X,Y)=\max_{1\le k\le n}\biggl\|\,\sum_{i=1}^{k}X_{i}-%
\sum_{i=1}^{k}Y_{i}\,\biggr\|
\]
would be so small as possible with large probability. Here \thinspace $\Vert
\cdot \Vert $ \thinspace is the Euclidean norm. It is clear that the vectors
even with the same distributions can be very far one from another.

In some sense this problem is one of the most important in probability
approximations because many well-known probability theorems can be
considered as consequences of results about strong approximation of
sequences of sums by corresponding Gaussian sequences. This is related to
the law of iterated logarithm, to several theorems about large deviations,
to the estimates for the rate of convergence of the Prokhorov distance in
the invariance principles (Prokhorov [19], Skorokhod [26], Borovkov [4]), as
well as to the Strassen-type approximations (Strassen [28], see, for
example, Cs\"{o}rg\H{o} and Hall [8]).

The rate for strong approximation in the one-dimensional invariance
principle was studied by many authors (see, e.g., Prokhorov~[19],
Skorokhod~[26], Borovkov~[4], Cs\"{o}rg\H{o} and R\'{e}v\'{e}sz~[6] and the
bibliography in~Cs\"{o}rg\H{o} and R\'{e}v\'{e}sz~[7], Cs\"{o}rg\H{o} and
Hall~[8], Shao~[20]). Skorokhod~[26] developed a method of construction of
close sequences of sequential sums of independent random variables on the
same probability space. For a long time the best rates of approximation were
obtained by this method, known now as the Skorokhod embedding. However,
Koml\'{o}s, Major and Tusn\'{a}dy (KMT)~[17] elaborated a new, more powerful
method of dyadic approximation. With the help of this method they obtained
optimal rates of Gaussian approximation for sequences of independent
identically distributed random variables.

We restrict ourselves on the most important case, where the summands have
finite exponential moments. Sakhanenko~[24] generalized and essentially
sharpened KMT results in the case of non-identically distributed random
variables. He considered the following class of one-dimensional
distributions:
\[
\mathcal{S}_{1}(\tau )=\left\{ \mathcal{L}(\xi ):\ \mathbf{E}\,\xi =0,\
\mathbf{E}\,|\xi |^{3}\exp \left( \tau ^{-1}|\xi |\right) \le \tau \,\mathbf{%
E}\,|\xi |^{2}\right\}
\]
(the distribution of a random vector $\xi $ will be denoted by $\mathcal{L}%
(\xi ))$. His main result is formulated as follows.\medskip

\textbf{Theorem 1}{\ \textrm{(Sakhanenko [24])\bf .}} \textit{Suppose that }$%
\tau >0$\textit{, and }$\xi _{1},\dots ,\xi _{n}$\textit{\ are independent
random variables with \thinspace }$\mathcal{L}(\xi _{j})\in \mathcal{S}%
_{1}(\tau )$\textit{, \thinspace }$j=1,\dots ,n$\textit{. \thinspace Then one can construct on a probability space
a sequence of independent random variables \thinspace }$X_{1},\dots ,X_{n}$\textit{\ \thinspace and a sequence of
independent Gaussian random variables \thinspace }$Y_{1},\dots ,Y_{n}$\textit{\ \thinspace so that
}$\mathcal{L}(X_{j})$ $=\mathcal{L}(\xi
_{j})$\textit{, \ }$\mathbf{E}\,Y_{j}=0$\textit{, \ }${\mathbf{E}\,Y_{j}^{2}=%
\mathbf{E}\,X_{j}^{2}}$\textit{, \ }$j=1,\dots ,n$\textit{, \ and }
\begin{equation}
\mathbf{E}\,\exp \left( {c\,\Delta (X,Y)}/{\tau }\right) {\le 1+B/\tau ,}
\label{expsah}
\end{equation}
\textit{where }$c$\textit{\ is an absolute constant and \ }$B^{2}=\mathbf{E}%
\,\xi _{1}^{2}+\dots +\mathbf{E}\,\xi _{n}^{2}$\textit{.}\medskip

KMT [17] supposed that $\xi ,\xi _{1},\dots ,\xi _{n}$ are identically
distributed and \thinspace $\mathbf{E}\,e^{\langle h,\xi \rangle }<\infty $,
\thinspace for \thinspace $h\in V$, \thinspace where $V\subset \mathbf{R}^{d}
$ \thinspace is some neighborhood of zero. The KMT (1975--76) result follows
from Theorem 1. It is easy to see that there exists $\tau (F)$ such that $F=%
\mathcal{L}(\xi _{j})\in \mathcal{S}_{1}(\tau (F))$. Applying the Chebyshev
inequality, we observe that (\ref{expsah}) imply that
\begin{equation}
\mathbf{P}\,\left( c_{1}\,\Delta (X,Y)/\tau (F)\geq x\right) \le \exp \left(
\log \left( 1+\sqrt{{n}\,\mathbf{E}\,\xi ^{2}}/\tau (F)\right) -x\right) {%
,\quad x>0.}  \label{KMT}
\end{equation}
Inequality (\ref{KMT}) provides more information than the original KMT
formulation which contains unspecified constants depending on $F$. In (\ref
{KMT}) the dependence of constants on the distribution $F$ is written out in
an explicit form. The quantity $\tau (F)$ can be easily calculated or
estimated for any concrete distribution $F$.

The first attempts to extend the KMT and Sakhanenko approximations to the
multidimensional case (see Berkes and Philipp~[3], Philipp~[18], Berger~[2],
Einmahl~[10, 11]) had a partial success only. Comparatively recently U.
Einmahl [12] obtained multidimensional analogs of KMT results which are
close to optimal. Zaitsev [33, 34] removed an unnecessary logarithmic factor
from the result of Einmahl~[12] and obtained multidimensional analogs of KMT
results (see Theorem 2 below). In Theorem 2 the random vectors are,
generally speaking, non-identically distributed. However, they have the same
identity covariance operator~\thinspace $\Bbb{I}$. \thinspace Therefore, the
problem of obtaining an adequate multidimensional generalization of the main
result of Sakhanenko~[24] remained open. This generalization is given in
Theorem~3 below.

\section{Main results} \label{section 2} \setzero

\vskip-5mm \hspace{5mm}

For formulations of results we need some notations. Let $\mathcal{A}%
_{d}(\tau )$, \thinspace $\tau \geq 0$, $d\in \mathbf{N}$, denote classes of
$d$-dimensional distributions, \thinspace introduced in Zaitsev [29], see as
well Zaitsev [33--35]. The class \thinspace $\mathcal{A}_{d}(\tau )$
\thinspace (with a fixed \thinspace $\tau \geq 0$) consists of $d$%
-dimensional distributions~$F$ \thinspace for which the function $\varphi
(z)=\varphi (F,z)=\log \int_{\mathbf{R}^{d}}e^{\left\langle z,x\right\rangle
}F\{dx\}$ $(\varphi (0)=0)$ is defined and analytic for \thinspace $\left\|
z\right\| \tau <1$, $z\in \mathbf{C}^{d}$, \thinspace and $\left|
d_{u}d_{v}^{2}\,\varphi (z)\right| \leq \Vert u\Vert \tau \,\left\langle
\Bbb{D}\,v,v\right\rangle $ for all $u,v\in \mathbf{R}^{d}$ and $\left\|
z\right\| \tau <1$, where \thinspace $\Bbb{D}=\mbox{cov}\,F$, the covariance
operator corresponding to $F$, and $d_{u}\varphi $ is the derivative
\thinspace of the function $\varphi $ in direction $u$.\medskip

\textbf{Theorem 2 }(Zaitsev, [33, 34])\textbf{\textrm{.}} \textit{Suppose
that \thinspace }$\tau \ge 1$\textit{, \thinspace }$\alpha >0$\textit{\
\thinspace and }$\xi _{1},\dots ,\xi _{n}$\textit{\ \thinspace are random
vectors with distributions \thinspace }${\mathcal{L}(\xi _{k})\in \mathcal{A}%
_{d}(\tau )}$\textit{, }$\mathbf{E}\,\xi _{k}=0$\textit{, }$\mbox{cov}\xi
_{k}=\Bbb{I}$\textit{, \thinspace \thinspace }$k=1,\dots ,n$\textit{.
\thinspace Then one can construct on a probability space a sequence of
independent random vectors }$X_{1},\dots ,X_{n}$\textit{\ and a sequence of
independent Gaussian random vectors }$\,Y_{1},\dots ,Y_{n}\,$\textit{\ so
that }
\[
\mathcal{L}(X_{k})=\mathcal{L}(\xi _{k}),\quad \mathbf{E}\,Y_{k}=0,\quad %
\mbox{cov}\,{\mathcal{L}}\left( Y_{k}\right) =\Bbb{I},\qquad k=1,\dots ,n,
\]
\textit{and }
\[
\mathbf{E}\,\exp \left( \frac{c_{1}(\alpha )\,\Delta (X,Y)}{\tau d^{7/2}\log
^{*}d}\right) \le \exp \left( c_{2}(\alpha )\,d^{9/4+\alpha }\log ^{*}\left(
n/\tau ^{2}\right) \right) ,
\]
\textit{where }$c_{1}(\alpha )$\textit{, }$c_{2}(\alpha )$\textit{\ are positive quantities depending
on~\thinspace }$\alpha $\textit{\ only and }$\log ^{*}b=\max \{ 1$, $\log b\} $\textit{, for }$b>0.$\medskip

\textbf{Corollary 1.} \textit{In the conditions of Theorem~}$2$\textit{\ for
all \thinspace }$x\ge 0$\textit{\ \thinspace the following inequality is
valid }
\[
\mathbf{P}\left\{ \,\Delta (X,Y)\ge \frac{c_{2}(\alpha )\,\tau
d^{23/4+\alpha }\log ^{*}d\,\log ^{*}\left( n/\tau ^{2}\right) }{%
c_{1}(\alpha )}+x\right\} \le \exp \left( -\frac{c_{1}(\alpha )\,x}{\tau
d^{7/2}\log ^{*}d}\right) .
\]
\medskip

It is easy to see that if \thinspace $V\subset \mathbf{R}^{d}$ \thinspace is
some neighborhood of zero and \thinspace $\mathbf{E}\,e^{\langle h,\xi
\rangle }<\infty $, \thinspace for \thinspace $h\in V$, \thinspace then $F=%
\mathcal{L}(\xi )\in \mathcal{A}_{d}(c(F){)}$. Below we list some simple and
useful properties of classes $\mathcal{A}_{d}(\tau )$ which are essential in
the proof of Theorem 2. Theorem 2 implies in one-dimensional case
Sakhanenko's Theorem 1 for identically distributed random variables with
finite exponential moments as well as the result of KMT [17].\medskip

\textbf{Corollary 2.} \textit{Suppose that a random vector~\thinspace }$\xi $%
\textit{\ \thinspace has finite exponential moments \thinspace }$\mathbf{E}%
\,e^{\langle h,\xi \rangle }$\textit{, \thinspace for \thinspace }$h\in V$%
\textit{, \thinspace where \thinspace }$V\subset \mathbf{R}^{d}$\textit{\
\thinspace is some neighborhood of zero. Then one can construct on a
probability space a sequence of independent random vectors \thinspace }$%
X_{1},X_{2},\dots $\textit{\ \thinspace and a sequence of independent
Gaussian random vectors \thinspace }$Y_{1},Y_{2},\dots $\textit{\ \thinspace
so that }
\[
\mathcal{L}(X_{k})=\mathcal{L}(\xi ),\quad \mathbf{E}\,Y_{k}=0,\quad %
\mbox{cov}Y_{k}=\mbox{cov}\xi ,\qquad k=1,2,\dots ,
\]
\textit{and }
\[
\sum_{k=1}^{n}X_{k}-\sum_{k=1}^{n}Y_{k}=O(\log n)\qquad
\mbox{a.s.\,.}
\]
\medskip

As it is noted in KMT [17], from the results of B\'{a}rtfai [1] that the
rate of approximation in Corollary~2 is the best possible for non-Gaussian
vectors~$\xi $. An analog of Corollary 2 was obtained by Einmahl [12] under
additional smoothness-type restrictions on the distribution $\mathcal{L}(\xi
)$. The following statement is a sharpening of Corollary 2.\medskip

\textbf{Corollary 3 }(Zaitsev [36])\textbf{.} \textit{Suppose that a random
vector~\thinspace }$\xi $\textit{\ \thinspace has the distribution such that
}$\mathcal{L}(\Bbb{D}^{-1/2}\xi )\in \mathcal{A}_{d}(\tau )$\textit{, where }%
$\Bbb{D}=\mbox{cov}\,\mathcal{L}(\xi )$ \textit{is a reversible operator.
Let }$\sigma ^{2}$, $\sigma >0$, \textit{be the maximal eigenvalue of }$\Bbb{%
D}$. \textit{Then for any }$\alpha >0$ \textit{there exists a construction
from Corollary }$2$ \textit{such that}
\begin{equation}
\mathbf{P}\left\{ \limsup_{n\rightarrow \infty }\frac{1}{\log n}\left\|
\sum_{k=1}^{n}X_{k}-\sum_{k=1}^{n}Y_{k}\right\| \leq c_{3}(\alpha )\,\sigma
\,\tau \,d^{23/4+\alpha }\log ^{*}d\right\} =1  \label{cor3}
\end{equation}
\textit{with }$c_{3}(\alpha )$\textit{\ depending on }$\alpha $\textit{\
only.}\medskip

In Theorems 2 and Corollary 3 we consider the case $\tau \ge 1$. The case of
small $\tau $ was investigated by G\"{o}tze and Zaitsev [16]. It is shown
that under additional smoothness-type restrictions on the distribution $%
\mathcal{L}(\xi )$ the expression in the right-hand side of the inequality
in (\ref{cor3}) can be arbitrarily small if the parameter $\tau $ is small
enough. It is clear that the statements of Theorem~2 and Corollary 3 becomes
stronger for small~\thinspace $\tau $. \thinspace In G\"{o}tze and Zaitsev
[16] one can find simple examples in which the sufficiently complicated
smoothness condition is satisfied. The approximation is better in the case
when summands have smooth distributions which are close to Gaussian ones
(see inequalities (\ref{pi0}) and (\ref{pi}) below).

The following Theorem 3 is a generalization of Theorem 2 to the case of
multivariate random variables. In one-dimensional situation, Theorem 3
implies Theorem~1.\textit{\medskip }

\textbf{Theorem 3}{\ \textrm{(Zaitsev [35])\bf .}} \textit{Suppose
that }$\alpha
>0$\textit{, }$\tau \ge 1$\textit{, and }$\xi _{1},\dots ,\xi _{n}$\textit{\
are independent random vectors with \thinspace }$\mathbf{E}\,\xi _{j}=0$%
\textit{, \thinspace }$j=1,\dots ,n$\textit{. \thinspace Assume that there
exists a strictly increasing sequence of non-negative integers \thinspace }$%
m_{0}=0$\textit{, }$m_{1},\dots ,m_{s}=n$\textit{\ \thinspace satisfying the
following conditions. Write }
\[
\zeta _{k}=\xi _{m_{k-1}+1}+\dots +\xi _{m_{k}},\qquad k=1,\dots ,s,
\]
\textit{and suppose that }$($\textit{for all \thinspace }$k=1,\dots ,s)$%
\textit{\ }${\mathcal{L}(\zeta _{k})\in \mathcal{A}_{d}(\tau )}$\textit{, }$%
\mbox{cov}\,\zeta _{k}=\Bbb{B}_{k}$\textit{\ and, for all }$u\in \mathbf{R}%
^{d}$\textit{, }
\begin{equation}
c_{4}\left\| u\right\| ^{2}\le \left\langle \Bbb{B}_{k}u,u\right\rangle \le
c_{5}\left\| u\right\| ^{2}  \label{c}
\end{equation}
\textit{with some constants \thinspace }$c_{4}$\textit{\ \thinspace and
\thinspace }$c_{5}$\textit{. Then one can construct on a probability space a
sequence of independent random vectors \thinspace }$X_{1},\dots ,X_{n}$%
\textit{\ \thinspace and a corresponding sequence of independent Gaussian
random vectors \thinspace }$Y_{1},\dots ,Y_{n}$\textit{\ \thinspace so that }%
$\mathcal{L}(X_{j})=\mathcal{L}(\xi _{j})$\textit{, }$\mathbf{E}\,Y_{j}=0$%
\textit{, }${\mbox{cov}\,\mathcal{L}(Y_{j})=\mbox{cov}\,\mathcal{L}(X_{j})}$%
\textit{, }$j=1,\dots ,n$\textit{, and }
\[
\mathbf{E}\,\exp \left( \frac{a_{1}\,\Delta (X,Y)}{\tau d^{9/2}\log ^{*}d}%
\right) \le \exp \left( a_{2}\,d^{3+\alpha }\,\log ^{*}(s/\tau ^{2})\right)
,
\]
\textit{where }$a_{1}$\textit{, }$a_{2}$\textit{\ are positive quantities
depending only on \thinspace }$\alpha ,\,c_{4},\,c_{5}$\textit{.}\medskip

\section{Properties of classes {\boldmath $\mathcal{A}_{d}(\tau )$}}

\label{section 3} \setzero\vskip-5mm \hspace{5mm }

Let us consider elementary properties of classes $\mathcal{A}_{d}(\tau )$
which are essentially used in the proof of Theorems 2 and 3, see Zaitsev
[29, 31, 33--35]. It is easy to see that \thinspace $\tau _{1}<\tau _{2}$
\thinspace implies \thinspace $\mathcal{A}_{d}\left( {\tau _{1}}\right)
\subset \mathcal{A}_{d}\left( {\tau _{2}}\right) $. \thinspace Moreover, the
class $\mathcal{A}_{d}\left( {\tau }\right) $ is closed with respect to
convolution: if~\thinspace $F_{1},F_{2}\in \mathcal{A}_{d}\left( {\tau }%
\right) $, then~$F_{1}F_{2}$\thinspace $=F_{1}*F_{2}\in \mathcal{A}%
_{d}\left( {\tau }\right) $. Products of measures are understood in the
convolution sense. Note that the condition ${\mathcal{L}(\zeta _{k})\in
\mathcal{A}_{d}({\tau })}$ in Theorem 3 is satisfied if $\mathcal{L}(\xi
_{j})\in \mathcal{A}_{d}(\tau )$, for \thinspace $j=1,\dots ,n$.

Let \thinspace $\tau \ge 0$, $F=\mathcal{L}(\xi )\in \mathcal{A}_{d}(\tau )$%
, $y\in \mathbf{R}^{m}$, and $\Bbb{A}:\mathbf{R}^{d}\to \mathbf{R}^{m}$ is a
linear operator. Then
\[
\mathcal{L}(\Bbb{A\xi +}y)\in \mathcal{A}_{m}\left( \left\| \Bbb{A}\right\|
\tau \right) ,\quad \mbox{where}\ \ \left\| \Bbb{A}\right\| =\sup_{x\in
\mathbf{R}^{d},\,\left\| x\right\| \le 1}\left\| \Bbb{A}x\right\| .
\]

Suppose that \thinspace $\tau \ge 0$, $F_{k}=\mathcal{L}\left( \xi
^{(k)}\right) \in \mathcal{A}_{d_{k}}(\tau )$, and the vectors $\xi ^{(k)}$,
$k=1,2$, are independent. Let $\xi \in \mathbf{R}^{d_{1}+d_{2}}$ be the
vector with the first \thinspace $d_{1}$~coordinates coinciding with those
of~\thinspace $\xi ^{(1)}$ and with the last \thinspace $d_{2}$~coordinates
coinciding with those of $\xi ^{(2)}$. Then $F=\mathcal{L}(\xi )\in \mathcal{%
A}_{{{d_{1}+d_{2}}}}({\tau )}$.

The classes \thinspace $\mathcal{A}_{d}(\tau )$ \thinspace are closely
connected with other naturally defined classes of multidimensional
distributions. From the definition of \thinspace $\mathcal{A}_{d}(\tau )$ it
follows that if \thinspace ${\mathcal{L}(\xi )\in \mathcal{A}_{d}(\tau )}$
then the vector~\thinspace $\xi $ \thinspace has finite exponential moments
\thinspace $\mathbf{E}\,e^{\langle h,\xi \rangle }<\infty $, \thinspace for
\thinspace \thinspace $h\in \mathbf{R}^{d}$, \thinspace $\left\| h\right\|
\tau <1$. \thinspace This leads to exponential estimates for the tails of
distributions.

The condition $\mathcal{L}(\xi )\in \mathcal{A}_{{1}}({\tau )}$ is
equivalent to Statulevi\v {c}ius' [27] conditions on the rate of increasing
of cumulants \thinspace $\gamma _{m}$ \thinspace of the random variable $\xi
$:
\[
|\gamma _{m}|\le \frac{1}{2}m!\,\tau ^{m-2}\gamma _{2},\qquad m=3,4,\dots .
\]
This equivalence means that if one of these conditions is satisfied with
parameter~\thinspace $\tau $, \thinspace then the second is valid with
parameter~\thinspace $c\tau $, where $c$ denotes an absolute constant.
\thinspace However, the condition \thinspace $\mathcal{L}(\xi )\in \mathcal{A%
}_{{d}}({\tau )}$ differs essentially from other multidimensional analogs of
Statulevi\v {c}ius' conditions, considered by Rudzkis~[23] and Saulis~[25].

Zaitsev [30] considered classes of distributions
\begin{eqnarray*}
\mathcal{B}_{d}(\tau ) &=&\Big\{F=\mathcal{L}(\xi ):\mathbf{E}\,\xi =0,\
\left| \mathbf{E}\,\left\langle \xi ,v\right\rangle ^{2}\left\langle \xi
,u\right\rangle ^{m-2}\right|  \\
&\le &\frac{1}{2}m!\,\tau ^{m-2}\left\| u\right\| ^{m-2}\,\mathbf{E}%
\,\left\langle \xi ,v\right\rangle ^{2}\ \ \mbox{for all}\ u,v\in \mathbf{R}%
^{d},\ m=3,4,\dots \Big\}
\end{eqnarray*}
satisfying multidimensional analogs of the Bernstein inequality condition.
Sakhanenko's condition $\mathcal{L}(\xi )\in \mathcal{S}_{{1}}({\tau )}$ is
equivalent to the condition $\mathcal{L}(\xi )\in \mathcal{B}_{{1}}({\tau )}$%
. Note that if $F\left\{ \left\{ x\in \mathbf{R}^{d}:\left\| x\right\| \le
\tau \right\} \right\} =1\ $then $F\in \mathcal{B}_{d}(\tau ).$

Let us formulate a relation between classes \thinspace $\mathcal{A}_{d}(\tau
)$ \thinspace and $\mathcal{B}_{d}(\tau )$. Denote by $\sigma ^{2}(F)$ the
maximal eigenvalue of the covariance operator of a~distribution~\thinspace $%
F $. Then

\medbreak

\noindent a) If $F=\mathcal{L}(\xi )\in \mathcal{B}_{d}(\tau )$, then $%
\sigma ^{2}(F)\le 12\,\tau ^{2}$, $\mathbf{E}\,\xi =0$ and ${F\in \mathcal{A}%
_{d}(c\tau )} $.

\noindent b) If ${F=\mathcal{L}(\xi )\in \mathcal{A}_{d}(\tau )}$, $\sigma
^{2}(F)\le \tau ^{2}$ and $\mathbf{E}\,\xi =0$, then ${F\in \mathcal{B}%
_{d}(c\tau )}$.

\medbreak

If \thinspace $F$ \thinspace is an infinitely divisible distributions with
spectral measure concentrated on the ball $\left\{ x\in \mathbf{R}%
^{d}:\left\| x\right\| \le \tau \right\} $ then ${F\in \mathcal{A}_{d}(c\tau
)}$, where $c$ is an absolute constant. It is obvious that the class $%
\mathcal{A}_{d}\left( {0}\right) $ coincides with the class of all $d$%
-dimensional Gaussian distributions. The following inequality was proved
in~Zaitsev [29] and can be considered as an estimate of stability of this
characterization:
\begin{equation}
\mbox{if}\ F\in \mathcal{A}_{d}(\tau ),\ \mbox{then}\ \pi \left( F,\,\Phi
(F)\right) \le c\,d^{2}\tau \,\log ^{*}(\tau ^{-1});  \label{pi0}
\end{equation}
where \thinspace $\pi (\cdot ,\cdot )$ \thinspace is the Prokhorov distance
and $\Phi (F)$ denotes the Gaussian distribution whose mean and covariance
operator coincide with those of~\thinspace $F$. The Prokhorov distance
between distributions $F,G$ may be defined by means of the formula
\[
\pi (F,G)=\inf \left\{ \lambda :\pi (F,G,\lambda )\leq \lambda \right\} ,
\]
where
\[
\pi (F,G,\lambda )=\sup_{X}\max \left\{ F\{X\}-G\{X^{\lambda
}\},G\{X\}-F\{X^{\lambda }\}\right\} ,\quad \lambda >0,
\]
and $X^{\lambda }=\{y\in \mathbf{R}^{d}:\inf\limits_{x\in X}\left\|
x-y\right\| <\lambda \}$ is the $\lambda $-neighborhood of the Borel set~$X$%
. Moreover, in~Zaitsev [29] it was established that
\begin{equation}
\pi (F,\Phi (F),\lambda )\le c\,d^{2}\exp \Big(-\frac{\lambda }{c\,d^{2}\tau
}\Big).  \label{pi}
\end{equation}
It~is very essential (and important) that the inequality~(\ref{pi}) is
proved for all~\thinspace ${\tau >0}$ \thinspace and for
arbitrary~\thinspace cov$F$, \thinspace in contrast to Theorems~2 and~3,
where \thinspace $\tau \ge 1$ \thinspace and covariance operators satisfy
condition (\ref{c}). \thinspace The question about the necessity of
condition~ (\ref{c}) in Theorems~2 and~3 remains open. In~Zaitsev [30]
inequalities (\ref{pi0}) and (\ref{pi}) were proved for convolutions of
distributions from $\mathcal{B}_{d}(\tau )$

By the Strassen--Dudley theorem (see Dudley [9]) coupled with inequality (%
\ref{pi}), one can construct on a probability space the random vectors
\thinspace $\xi $ \thinspace and \thinspace $\eta $ \thinspace with
\thinspace $\mathcal{L}(\xi )=F$ \thinspace and \thinspace $\mathcal{L}(\eta
)=\Phi (F)$ \thinspace so that
\begin{equation}
\mathbf{P}\left\{ \Vert \xi -\eta \Vert >\lambda \right\} \le c\,d^{2}\exp %
\Big(-\frac{\lambda }{c\,d^{2}\tau }\Big).  \label{pi1}
\end{equation}
For convolutions of bounded measures, this fact was used by Rio~[21],
Einmahl and Mason [13], Bovier and Mason [5], Gentz and L\"{o}we [15],
Einmahl and Kuelbs [14].

The scheme of the proof of Theorems~2 and 3 is very close to that of the
main results of Sakhanenko [24] and Einmahl [12]. We suppose that the
Gaussian vectors \thinspace $Y_{1},\dots ,Y_{n}$, \thinspace $n=2^{N}$,
\thinspace are already constructed and construct the independent vectors
which are bounded with probability one, have sufficiently smooth
distributions and the same moments of the first, second and third orders as
the needed independent random vectors~\thinspace $X_{1},\dots ,X_{n}$.
\thinspace For the construction we use the dyadic scheme proposed by KMT
[17]. Firstly we construct the sum of $2^{N}$ summands using the Rosenblatt
[22] quantile transform for conditional distributions (see Einmahl [12]).
Then we construct blocks of $2^{N-1},2^{N-2},\ldots ,1$ summands. The rate
of approximation is estimated using the fact that, for smooth summands
distributions, the corresponding conditional distribution are smooth and
close to Gaussian ones. Then we construct the vectors~ $X_{1},\dots ,X_{n}$%
\thinspace \thinspace in several steps. After each step the number
of~\thinspace $X_{k}$ \thinspace which are not constructed becomes smaller
in \thinspace $2^{p}$ ~\thinspace times, where \thinspace $p$ ~\thinspace is
a suitably chosen positive integer. In each step we begin with already
constructed vectors which are bounded with probability one and have
sufficiently smooth distributions and the needed moments up to the third
order. Then we construct the vectors such that, in each block of \thinspace $%
2^{p}$ ~\thinspace summands, only the first vector has the initial bounded
smooth distribution. The rest \thinspace $2^{p}-1$ ~\thinspace vectors have
the needed distributions ~\thinspace $\mathcal{L}(\xi _{k})$. \thinspace
These \thinspace $2^{p}-1$ ~\thinspace vectors from each block will be
chosen as $X_{k}$ \thinspace and will be not involved in the next steps of
the procedure. The coincidence of third moments will allow us to use more
precise estimates of the closeness of quantiles of conditional distributions
contained in Zaitsev [32]. In the estimation of closeness of random vectors
in the steps of the procedure described above, we use essentially properties
of classes $\mathcal{A}_{d}(\tau )$.

\section{Infinitely divisible approximation}

\label{section 4} \setzero\vskip-5mm \hspace{5mm}

Let us finally mention a result about strong approximation of sums of
independent random vectors by infinitely divisible distributions. Theorem 4
below follows from the main result of Zaitsev [32] coupled with the
Strassen--Dudley theorem. Inequality (\ref{infdiv}) can be considered as a
generalization of inequality (\ref{pi1}) to convolution of distribution with
unbounded supports.\medskip

\textbf{Theorem 4}\textrm{.} \textit{Let }$d$\textit{-dimensional
probability distributions }$F_{i}$\textit{, }$i=1,\dots ,n$\textit{, be
represented as mixtures of }$d$\textit{-dimensional probability
distributions \thinspace }$U_{i}$\textit{\ \thinspace and \thinspace }$V_{i}$%
\textit{: }
\[
F_{i}=(1-p_{i})U_{i}+p_{i}V_{i},
\]
\textit{where }
\[
0\le p_{i}\le 1,\qquad \int x\,U_{i}\{dx\}=0,\qquad U_{i}\left\{ \left\{
x\in \mathbf{R}^{d}:\left\| x\right\| \le \tau \right\} \right\} =1,
\]
\textit{and }$V_{i}$\textit{\ are arbitrary distributions. Then for any
fixed }$\lambda >0$\textit{\ one can construct on the same probability space
the random vectors \thinspace }$\xi $\textit{\ \thinspace and \thinspace }$%
\eta $\textit{\ \thinspace so that }
\begin{equation}
\mathbf{P}\left\{ \Vert \xi -\eta \Vert >\lambda \right\} \le c(d)\,\left(
\max_{1\le i\le n}p_{i}+\exp \left( -\frac{\lambda }{c(d)\tau }\right)
\right) +\sum_{i=1}^{n}p_{i}^{2}  \label{infdiv}
\end{equation}
\textit{and }
\[
\mathcal{L}(\xi )=\prod_{i=1}^{n}F_{i},\qquad \mathcal{L}(\xi
)=\prod_{i=1}^{n}\mbox{e}(F_{i}),
\]
\textit{where }$c(d)$\textit{\ depends on only and \thinspace }$\mbox{e}%
(F_{i})$\textit{\ \thinspace denotes the compound Poisson infinitely
divisible distribution with characteristic function \thinspace }$\exp (%
\widehat{F}_{i}(t)-1)$\textit{, where }$\widehat{F}_{i}(t)=\int e^{itx}\,F_{i}\{dx\}$\textit{. If the
distributions }$V_{i}$\textit{\ are identical, the term \thinspace }$\sum_{i=1}^{n}p_{i}^{2}$\textit{\ \thinspace
in }(\ref{infdiv})\textit{\ can be omitted}.

\label{lastpage}

\end{document}